\documentclass[leqno,12pt]{article}
\usepackage{amsmath,amssymb,amsfonts,mathrsfs,exscale,colonequals,a4wide,color}
\usepackage[curve,matrix,arrow,cmtip]{xy}
\usepackage{verbatim}
\NoComputerModernTips
\newdir^{ (}{{}*!/-3pt/\dir^{(}}   % Usage: \ar@{^{ (}->}[r]^{j}
\newdir_{ (}{{}*!/-3pt/\dir_{(}}   % Usage: \ar@{_{ (}->}[r]^{j}    
\entrymodifiers={+!!<0pt,\fontdimen22\textfont2>}

\newcommand{\defeq}{\colonequals}
\def\Tildem{{\widetilde{{\vrule width0pt height 6.8pt\smash{\widetilde m}}}}}
\def\tildem{{\widetilde{m}}}

\def\Aut{\mathop{\rm Aut}\nolimits}

\def\Hom{\mathop{\rm Hom}\nolimits}

\def\Gal{\mathop{\rm Gal}\nolimits}
\def\End{\mathop{\rm End}\nolimits}

\def\Mat{\mathop{\rm Mat}\nolimits}

\def\Quot{\mathop{\rm Quot}\nolimits}
\def\Div{\mathop{\rm Div}\nolimits}

\def\GL{\mathop{\rm GL}\nolimits}
\def\SL{\mathop{\rm SL}\nolimits}

\def\ad{{\rm ad}}

\def\ad{{\rm ad}}

\def\sep{{\rm sep}}

\def\opp{{\rm opp}}

\let\phi\varphi
\let\epsilon\varepsilon
\let\setminus\smallsetminus

\let\le\leqslant
\let\ge\geqslant

\newtheorem{Thm}{Theorem}[section]
\newtheorem{Prop}[Thm]{Proposition}
\newtheorem{Lem}[Thm]{Lemma}

\newtheorem{Cor}[Thm]{Corollary}

\newtheorem{Def}[Thm]{Definition}

\numberwithin{equation}{section}
\def\UseTheoremCounterForNextEquation{\setcounter{equation}{\value{Thm}}\addtocounter{Thm}{1}}
\def\qed{{\hskip0pt\unskip\unskip\nobreak\hfil\penalty50
          \hskip1em\hbox{}\nobreak\hfil
           {$\square$}
          \parfillskip=0pt\finalhyphendemerits=0
          \par}\medskip}

\newenvironment{Proof}
               {\noindent{\bf Proof.}\ }
               {\qed}

               {\noindent{\bf Proof of #1.}\ }
               {\qed}

\newcommand{\BF}{{\mathbb{F}}}
\newcommand{\BG}{{\mathbb{G}}}

\newcommand{\Fa}{{\mathfrak{a}}}
\newcommand{\Fb}{{\mathfrak{b}}}
\newcommand{\Fc}{{\mathfrak{c}}}
\newcommand{\Fd}{{\mathfrak{d}}}

\newcommand{\Fm}{{\mathfrak{m}}}

\newcommand{\Fp}{{\mathfrak{p}}}
\newcommand{\Fq}{{\mathfrak{q}}}

\let\longto\longrightarrow

\let\into\hookrightarrow

\let\onto\twoheadrightarrow
\definecolor{Pink}{rgb}{1,0.2,0.7}

\begin{document}

%%%%%%%%%%%%%%%%%%%%%%%%%%%%%%%%%%%%%%%%%%%%%%%%%%%%%%%%%%%%%%%%%%%%%%%%%%%%%%%%%%%%%%%%%%%%%

\title{Kummer Theory for Drinfeld Modules}

\author{%
\renewcommand{\thefootnote}{\arabic{footnote}}%
Richard Pink%
\footnote{Dept. of Mathematics, ETH Z\"urich, CH-8092 Z\"urich, Switzerland,
\tt pink@math.ethz.ch} 
}

\date{\today}
\maketitle

\abstract{Let $\phi$ be a Drinfeld $A$-module of characteristic $\Fp_0$ over a finitely generated field~$K$. Previous articles determined the image of the absolute Galois group of $K$ up to commensurability in its action on all prime-to-$\Fp_0$ torsion points of~$\phi$, or equivalently, on the prime-to-$\Fp_0$ adelic Tate module of~$\phi$. In this article we consider in addition a finitely generated torsion free $A$-submodule $M$ of $K$ for the action of $A$ through~$\phi$. 
We determine the image of the absolute Galois group of $K$ up to commensurability in its action on the prime-to-$\Fp_0$ division hull of~$M$, or equivalently, on the extended prime-to-$\Fp_0$ adelic Tate module associated to $\phi$ and~$M$.}
% \abstract{For any Drinfeld $A$-module over a finitely generated field~$K$ of characteristic $\Fp_0$ and any finitely generated torsion free $A$-submodule $M$ of $K$ we determine the image of the absolute Galois group $\Gal(K^\sep/K)$ on the prime-to-$\Fp_0$ division hull of~$M$ up to commensurability.}

\bigskip\bigskip
{\parskip=0pt
\tableofcontents
}

\newpage

%%%%%%%%%%%%%%%%%%%%%%%%%%%%%%%%%%%%%%%%%%%%%%%%%%%%%%%%%%%%%%%%%%%%%%%%%%%%%%%%%%%%%%%%%%%%%

% \newpage
% \addtocounter{section}{-1}
\section{Introduction}
\label{Intro}

Let $F$ be a finitely generated field of transcendence degree $1$ over the prime field $\BF_p$ of characteristic $p>0$. Let $A$ be the ring of elements of $F$ which are regular outside a fixed place $\infty$ of~$F$. Let $K$ be another field that is finitely generated over~$\BF_p$, and let $K^\sep$ be a separable closure of~$K$. Write $\End(\BG_{a,K}) = K[\tau]$ with $\tau(x) = x^p$. Let $\phi\colon A\to K[\tau]$, $a\mapsto\phi_a$ be a Drinfeld $A$-module of rank $r\ge1$ and characteristic~$\Fp_0$. Then either $\Fp_0$ is the zero ideal of~$A$ and $\phi$ is said to have generic characteristic; or $\Fp_0$ is a maximal ideal of~$A$ and $\phi$ is said to have special characteristic.
% Note that $\Fp_0$ is the zero ideal of $A$ if $\phi$ has generic characteristic; otherwise $\phi$ has special characteristic and $\Fp_0$ is a maximal ideal of~$A$.
 
For brevity we call any maximal ideal of $A$ a prime of~$A$. For any prime $\Fp\not=\Fp_0$ of $A$ the $\Fp$-adic Tate module $T_\Fp(\phi)$ is a free module of rank $r$ over the completion~$A_\Fp$, endowed with a continuous action of the Galois group $\Gal(K^\sep/K)$. The prime-to-$\Fp_0$ adelic Tate module $T_\ad(\phi) = \prod_{\Fp\not=\Fp_0}T_\Fp(\phi)$ is then a free module of rank $r$ over $A_\ad = \prod_{\Fp\not=\Fp_0} A_\Fp$ carrying a natural action of Galois. This action corresponds to a continuous homomorphism \UseTheoremCounterForNextEquation
\begin{equation}\label{GamadDef}
\Gal(K^\sep/K) \ \to\ \Aut_{A_\ad}(T_\ad(\phi)) \ \cong\ \GL_r(A_\ad).
\end{equation}
Its image $\Gamma_\ad$ was determined up to commensurability by Pink-R\"utsche \cite{PinkRuetsche2} and Devic-Pink \cite{DevicPink}; for special cases see  Theorems \ref{MainIntro} and \ref{GammaAdelicOpen} below.

\medskip
Let $M\subset K$ be a finitely generated torsion free $A$-submodule of rank $d$ for the action of $A$ through~$\phi$. Then there is an associated prime-to-$\Fp_0$ adelic Tate module $T_\ad(\phi,M)$,
% (see (\ref{TadMDef}) below), 
which is a free module of rank $r+d$ over~$A_\ad$ carrying a natural continuous action of $\Gal(K^\sep/K)$. This module lies in a natural Galois equivariant short exact sequence
\UseTheoremCounterForNextEquation
\begin{equation}\label{TadM1}
\xymatrix{ 0 \ar[r] & T_\ad(\phi) \ar[r] & T_\ad(\phi,M) \ar[r] & M\otimes_AA_\ad \ar[r] & 0\rlap{.}\\}
\end{equation}
Define $\Gamma_{\ad,M}$ as the image of the continuous homomorphism
\UseTheoremCounterForNextEquation
\begin{equation}\label{GamadMDef}
\Gal(K^\sep/K) \ \to\ \Aut_{A_\ad}(T_\ad(\phi,M)) \ \cong\ \GL_{r+d}(A_\ad).
\end{equation}
Then the restriction to $T_\ad(\phi)$ induces a surjective homomorphism $\Gamma_{\ad,M}\onto\Gamma_\ad$, whose kernel we denote by $\Delta_{\ad,M}$. Since the action on  $M\otimes_AA_\ad$ is trivial, there is a natural inclusion
\UseTheoremCounterForNextEquation
\begin{equation}\label{DeltaAdMInclus1}
\Delta_{\ad,M} \into \Hom_A(M,T_\ad(\phi)).
\end{equation}
Any splitting of the sequence (\ref{TadM1}) induces an inclusion into the semidirect product
\UseTheoremCounterForNextEquation
\begin{equation}\label{GammaAdMInclus}
\Gamma_{\ad,M} \into \Gamma_\ad \ltimes \Hom_A(M,T_\ad(\phi)).
\end{equation}
The aim of this article is to describe these subgroups up to commensurability. 

\medskip
In general the shape of these Galois groups is affected by the endomorphisms
of $\phi$ over $K^\sep$, and in special characteristic also by the endomorphisms of the restrictions of $\phi$ to all subrings of~$A$. Any general results therefore involve further definitions and notation. In this introduction we avoid these and mention only a special case; the general case is addressed by Theorems \ref{MainPrimitive} and \ref{MainCentral} and \ref{MainGeneralSpec}.
Parts (a) and (b) of the following result are due to Pink-R\"utsche \cite[Thm.$\;$0.1]{PinkRuetsche2}, resp.\ Devic-Pink \cite[Thm.$\;$1.1]{DevicPink}, and part (c) is a special case of Theorem \ref{MainPrimitive} below:

\begin{Thm}\label{MainIntro}
Assume that $\End_{K^\sep}(\phi) = A$, and in special characteristic also that $\End_{K^\sep}(\phi|B) = \nobreak A$ for every integrally closed infinite subring $B\subset A$. 
\begin{enumerate}
\item[(a)] If $\phi$ has generic characteristic, then $\Gamma_\ad$ is open in $\GL_r(A_\ad)$.
\item[(b)] If $\phi$ has special characteristic, then $\Gamma_\ad$ is commensurable with $\overline{\langle a_0\rangle} \cdot \SL_r(A_\ad)$ for some central element $a_0\in A$ that generates a positive power of~$\Fp_0$.
\item[(c)] 
The inclusions (\ref{DeltaAdMInclus1}) and (\ref{GammaAdMInclus}) are both open.
\end{enumerate}
\end{Thm}

The method used to prove Theorem \ref{MainIntro} (c) and its generalizations is an adaptation of the Kummer theory for semiabelian varieties from Ribet \cite{RibetKummer} and predecessors. The main ingredients are the above mentioned descriptions of~$\Gamma_\ad$ and Poonen's tameness result \cite{PoonenMordellWeil} concerning the structure of $K$ as an $A$-module via~$\phi$. 
A standard procedure would be to first prove corresponding results for $\Fp$-division points for almost all primes $\Fp\not=\Fp_0$ of~$A$, and for $\Fp$-power division points for all $\Fp\not=\Fp_0$, and then to combine these individual results by taking products, as in \cite{RibetKummer}, \cite{ChiLi2001}, \cite{Li2001}, 
\cite{PinkRuetsche2}, \cite{DevicPink}, \cite{Haeberli}. Instead, we have found a shorter way by doing everything adelically from the start. 
The core of the argument is the proof of Lemma \ref{MainPrim1}. Therein we avoid the explicit use of group cohomology by trivializing an implicit $1$-cocycle with the help of a suitable central element of~$\Gamma_\ad$. 
On first reading the readers may want to restrict their attention to the case of Theorem \ref{MainIntro}, which requires only Section \ref{Glob}, a little from Section \ref{Glob1b}, and Section \ref{PrimCase} with simplifications, avoiding Sections \ref{RedStep} and \ref{GenCase} entirely. Some of this was worked out in the Master's thesis of H\"aberli \cite{Haeberli}. Our results generalize those of Chi-Li \cite{ChiLi2001} and Li \cite{Li2001}.

\medskip
The notation and the assumptions of this introduction remain in force throughout the article.

\section{Extended Tate modules}
\label{Glob}

Following the usual convention in commutative algebra we let $A_{(\Fp_0)}\subset F$ denote the localization of $A$ at~$\Fp_0$; this is equal to $F$ if and only if $\phi$ has generic characteristic. Observe that there is a natural isomorphism of $A$-modules
\UseTheoremCounterForNextEquation
\begin{equation}\label{ATorDecomp1}
A_{(\Fp_0)}/A\ \cong\ \bigoplus_{\Fp\not=\Fp_0} F_\Fp/A_\Fp,
\end{equation}
where the product is extended over all maximal ideals $\Fp\not=\Fp_0$ of~$A$ and where $F_\Fp$ and $A_\Fp$ denote the corresponding completions of $F$ and~$A$. This induces a natural isomorphism for the prime-to-$\Fp_0$ adelic completion of $A$:
\UseTheoremCounterForNextEquation
\begin{equation}\label{ATorDecomp2}
A_\ad\ \defeq\ \End_A(A_{(\Fp_0)}/A)
\ \cong\ \prod_{\Fp\not=\Fp_0} A_\Fp.
\end{equation}
As a consequence, for any torsion $A$-module $X$ that is isomorphic to $(A_{(\Fp_0)}/A)^{\oplus n}$ for some integer~$n$, the construction
\UseTheoremCounterForNextEquation
\begin{equation}\label{TDef}
T(X) := \Hom_A(A_{(\Fp_0)}/A,X)
\end{equation} 
yields a free $A_\ad$-module of rank~$n$. Reciprocally $T(X)$ determines $X$ completely up to a natural isomorphism $X\cong T(X)\otimes_{A_\ad}(A_{(\Fp_0)}/A)$. Thus any $A$-linear group action on $X$ determines and is determined by the corresponding $A_\ad$-linear group action on~$T(X)$. Moreover, with
\UseTheoremCounterForNextEquation
\begin{equation}\label{TpDef}
T_\Fp(X) := \Hom_A(F_\Fp/A_\Fp,X)
\end{equation}
the decompositions (\ref{ATorDecomp1}) and (\ref{ATorDecomp2}) induce a decomposition 
\UseTheoremCounterForNextEquation
\begin{equation}\label{TDecomp}
T(X) \cong \prod_{\Fp\not=\Fp_0} T_\Fp(X)
\end{equation}
This will give a concise way of defining the $\Fp$-adic and adelic Tate modules associated to the given Drinfeld module~$\phi$.

%%%%%%%%%%%%%%%%%%%%%%%%%%%%%%%%%%%%%%%%%%%%%%%

\medskip
We view $K^\sep$ as an $A$-module with respect to the action $A\times K^\sep\to K^\sep$, $(a,x)\mapsto\phi_a(x)$ and are interested in certain submodules. One particular submodule is~$K$. Let $M$ be a finitely generated torsion free $A$-submodule of rank $d\ge0$ contained in~$K$. Then the \emph{prime-to-$\Fp_0$ division hull of $M$ in $K^\sep$} is the $A$-submodule
\UseTheoremCounterForNextEquation
\begin{equation}\label{DivHullPrimeDef} 
\Div_{K^\sep}^{(\Fp_0)}(M) \ \defeq\ 
\bigl\{x\in K^\sep\,\bigm|\,\exists a\in A\setminus\Fp_0: \phi_a(x)\in M\bigr\}.
\end{equation}
Let $\Div_K^{(\Fp_0)}(M)$ denote the intersection of $\Div_{K^\sep}^{(\Fp_0)}(M)$ with~$K$.
For later use we recall the following result proved by Poonen \cite[Lemma$\,$5]{PoonenMordellWeil} when $K$ is a global field and $\phi$ has generic characteristic, and by Wang \cite{WangMordellWeil} in general:

\begin{Thm}\label{lettThm}
% If $K$ is finitely generated over~$\BF_p$, then 
$[\Div_K^{(\Fp_0)}(M):M]$ is finite.
\end{Thm}

%%%%%%%%%%%%%%%%%%%%%%%%%%%%%%%%%%%%%%%%%%%%%%%

As a special case of the above, the prime-to-$\Fp_0$ division hull of the zero module $\Div_{K^\sep}^{(\Fp_0)}(\{0\})$ is the module of all prime-to-$\Fp_0$ torsion points of $\phi$ in~$K^\sep$.

\begin{Prop}\label{DivExSeq}
There is a natural short exact sequence of $A$-modules
$$\xymatrix{ 
0 \ar[r] 
& \Div_{K^\sep}^{(\Fp_0)}(\{0\}) \ar[r]
& \Div_{K^\sep}^{(\Fp_0)}(M) \ar[r]
& M\otimes_A A_{(\Fp_0)} \ar[r] & 1, \\}$$
where the map on the right hand side is described by $x\mapsto\phi_a(x)\otimes\frac{1}{a}$ for any $a\in A\setminus\Fp_0$ satisfying $\phi_a(x)\in M$. 
\end{Prop}

\begin{Proof}
The map is well-defined, because for any $x\in \Div_{K^\sep}^{(\Fp_0)}(M)$ and any $a$, $b\in A\setminus\Fp_0$ satisfying $\phi_a(x)$, $\phi_b(x)\in M$ we have
$$\xymatrix@R-10pt{ 
\phi_a(x)\otimes\frac{1}{a} \ =\ \phi_a(x)\otimes\frac{b}{ab} \ =\ 
\phi_b(\phi_a(x))\otimes\frac{1}{ab} \ar@<70pt>@{}[d]|\Vert \\
\phi_b(x)\otimes\frac{1}{b} \ =\ \phi_b(x)\otimes\frac{a}{ab} \ =\ 
\phi_a(\phi_b(x))\otimes\frac{1}{ab} \rlap{.} \\}$$
The map is clearly $A$-linear, and one easily shows that its kernel is $\Div_{K^\sep}^{(\Fp_0)}(\{0\})$; hence the sequence is left exact. The exactness on the right results from the fact that for any $a\in A\setminus\Fp_0$ and any $m\in M$, the polynomial equation $\phi_a(x)=m$ is separable of positive degree and therefore has a solution in~$K^\sep$.
\end{Proof}

Dividing by~$M$, the exact sequence from Proposition \ref{DivExSeq} yields a natural short exact sequence of $A$-modules
\UseTheoremCounterForNextEquation
\begin{equation}\label{DivExDiag1}
\vcenter{\xymatrix{
0 \ar[r] 
& \Div_{K^\sep}^{(\Fp_0)}(\{0\}) \ar[r] 
& \Div_{K^\sep}^{(\Fp_0)}(M)/M \ar[r] 
& M\!\otimes_A\!(A_{(\Fp_0)}/A) \ar[r] & 0\rlap{.}\\}}
\end{equation}
By the general theory of Drinfeld modules the module on the left is isomorphic to $(A_{(\Fp_0)}/A)^{\oplus r}$, where $r$ is the rank of~$\phi$. Using the functor $T$ from (\ref{TDef}), the \emph{prime-to-$\Fp_0$ adelic Tate module of $\phi$} can be described canonically as
\UseTheoremCounterForNextEquation
\begin{equation}\label{TadDef}
T_\ad(\phi)\ \defeq\ T\bigl(\Div_{K^\sep}^{(\Fp_0)}(\{0\})\bigr)
\end{equation}
and is a free $A_\ad$-module of rank~$r$. 
Since $M$ is a projective $A$-module of rank~$d$, the module on the right of (\ref{DivExDiag1}) is isomorphic to $(A_{(\Fp_0)}/A)^{\oplus d}$, and together it follows that the module in the middle is isomorphic to $(A_{(\Fp_0)}/A)^{\oplus(r+d)}$. The \emph{extended prime-to-$\Fp_0$ adelic Tate module of $\phi$ and $M$}
\UseTheoremCounterForNextEquation
\begin{equation}\label{TadMDef}
T_\ad(\phi,M)\ \defeq\ T\bigl(\Div_{K^\sep}^{(\Fp_0)}(M)/M\bigr)
\end{equation}
is therefore a free $A_\ad$-module of rank~$r+d$. Moreover, the exact sequence (\ref{DivExDiag1}) yields a natural short exact sequence of $A_\ad$-modules
\UseTheoremCounterForNextEquation
\begin{equation}\label{TadM2}
\xymatrix@C+10pt{ 0 \ar[r] & T_\ad(\phi) \ar[r] & T_\ad(\phi,M) \ar[r] & M\otimes_AA_\ad \ar[r] & 0.\\}
\end{equation}
All this decomposes uniquely as $T_\ad(\phi) = \prod_{\Fp\not=\Fp_0} T_\Fp(\phi)$ etc.\ as in (\ref{TDecomp}).

%%%%%%%%%%%%%%%%%%%%%%%%%%%%%%%%%%%%%%%%%%%%%%%

\medskip
By construction there is a natural continuous action of the Galois group $\Gal(K^\sep/K)$ on all modules and arrows in Proposition \ref{DivExSeq} and in (\ref{DivExDiag1}). This induces a continuous action on the short exact sequence (\ref{TadM2}), which in turn determines the former two by the following fact:

\begin{Prop}\label{MbyMAct}
The action of $\Gal(K^\sep/K)$ on $\Div_{K^\sep}^{(\Fp_0)}(M)$ is completely determined by the action on $T_\ad(\phi,M)$.
\end{Prop}

\begin{Proof}
For any $\sigma\in\Gal(K^\sep/K)$ the endomorphism $x\mapsto\sigma(x)-x$ of $\Div_{K^\sep}^{(\Fp_0)}(M)$ is trivial on $M$, because that module is contained in~$K$. Also, the image of this endomorphism is contained in $\Div_{K^\sep}^{(\Fp_0)}(\{0\})$, because for any $a\in A\setminus\Fp_0$ with $\phi_a(x)\in M$ we have $\phi_a(\sigma(x)-x) = \sigma(\phi_a(x))-\phi_a(x) \allowbreak =\nobreak 0$. Thus the endomorphism factors through a homomorphism $\Div_{K^\sep}^{(\Fp_0)}(M)/M \longto \Div_{K^\sep}^{(\Fp_0)}(\{0\})$. But by (\ref{DivExDiag1}) the latter homomorphism is determined completely by the action of $\sigma$ on $\Div_{K^\sep}^{(\Fp_0)}(M)/M$, and thus by the action of $\sigma$ on $T_\ad(\phi,M)$, as desired.
\end{Proof}

%%%%%%%%%%%%%%%%%%%%%%%%%%%%%%%%%%%%%%%%%%%%%%%

\medskip
Let $\Gamma_\ad$ and $\Gamma_{\ad,M}$ denote the images of $\Gal(K^\sep/K)$ acting on $T_\ad(\phi)$ and $T_\ad(\phi,M)$, as in (\ref{GamadDef}) and (\ref{GamadMDef}). Restricting to $T_\ad(\phi)$ induces a surjective homomorphism $\Gamma_{\ad,M}\onto\Gamma_\ad$, and we define $\Delta_{\ad,M}$ by the short exact sequence
\UseTheoremCounterForNextEquation
\begin{equation}\label{DelGamGam}
\xymatrix{ 1 \ar[r] & \Delta_{\ad,M} \ar[r] & \Gamma_{\ad,M} \ar[r] & \Gamma_{\ad} \ar[r] & 1\rlap{.}\\}
\end{equation}
For any $m\in M$ take an element $t\in T_\ad(\phi,M)$ with image $m\otimes1$ in $M\otimes_AA_\ad$. Since any $\delta\in\Delta_{\ad,M}$ acts trivially on $T_\ad(\phi)$, the difference $\delta(t)-t$ depends only on $\delta$ and~$m$. Since $\delta$ also acts trivially on $M\otimes_AA_\ad$, the difference lies in $T_\ad(\phi)$ and therefore defines a map
\UseTheoremCounterForNextEquation
\begin{equation}\label{PairingDef}
\Delta_{\ad,M} \times M \longrightarrow T_\ad(\phi),\ 
(\delta,m) \mapsto \langle\delta,m\rangle := \delta(t)-t.
\end{equation}
By direct calculation this map is additive in $\delta$ and $A$-linear in~$m$. By the construction of $\Delta_{\ad,M}$ the adjoint of the pairing (\ref{PairingDef}) is therefore a natural inclusion, already mentioned in (\ref{DeltaAdMInclus1}):
\UseTheoremCounterForNextEquation
\begin{equation}\label{DeltaAdMInclus2}
\Delta_{\ad,M} \ \into\ 
% \Hom_{A_\ad}(M\otimes_AA_\ad,T_\ad(\phi)) \ \cong\ 
\Hom_A(M,T_\ad(\phi)).
\end{equation}

%%%%%%%%%%%%%%%%%%%%%%%%%%%%%%%%%%%%%%%%%%%%%%%

Let $R := \End_K(\phi)$ denote the endomorphism ring of $\phi$ over~$K$. This is an $A$-order in a finite dimensional division algebra over~$F$. It acts naturally on $K$ and $K^\sep$ and therefore on $\Div_{K^\sep}^{(\Fp_0)}(\{0\})$ and $T_\ad(\phi)$, turning the latter two into modules over $R_\ad \defeq R\otimes_AA_\ad$. As this action commutes with the action of~$\Gamma_\ad$, it leads to an inclusion
\UseTheoremCounterForNextEquation
\begin{equation}\label{GammaEnd}
\Gamma_\ad\ \subset\ \Aut_{R_\ad}(T_\ad(\phi)).
\end{equation}
The decomposition (\ref{ATorDecomp2}) induces a decomposition $R_\ad = \prod_{\Fp\not=\Fp_0} R_\Fp$ where $R_\Fp \defeq R\otimes_AA_\Fp$ acts naturally on $T_\Fp(\phi)$.

If $M$ is an $R$-submodule of~$K$, then $R$ and hence $R_\ad$ also act on $\Div_{K^\sep}^{(\Fp_0)}(M)$ and $T_\ad(\phi,M)$, and these actions commute with the action of~$\Gamma_{\ad,M}$. The inclusion (\ref{DeltaAdMInclus2}) then factors through an inclusion
\UseTheoremCounterForNextEquation
\begin{equation}\label{DeltaAdMInclus2D}
\Delta_{\ad,M} \ \into\ \Hom_R(M,T_\ad(\phi)).
\end{equation}
Moreover, any $R$-equivariant splitting of the sequence (\ref{TadM2}) then induces an embedding into the semidirect product
\UseTheoremCounterForNextEquation
\begin{equation}\label{GammaAdMInclus2}
\Gamma_{\ad,M} \into \Gamma_\ad \ltimes \Hom_R(M,T_\ad(\phi)).
\end{equation}

%%%%%%%%%%%%%%%%%%%%%%%%%%%%%%%%%%%%%%%%%%%%%%%%%%%%%%%%%%%%%%%%%%%%%%%%%%%%%%%%%%%%%%%%%%%%%

\section{Reduction steps}
\label{RedStep}

For use in Section \ref{GenCase} we now discuss the behavior of extended Tate modules and their associated Galois groups under isogenies and under restriction of $\phi$ to subrings. 
% We keep the notations of Section~\ref{Glob}.

%%%%%%%%%%%%%%%%%%%%%%%%%%%%%%%%%%%%%%%%%%%%%%%
\medskip
First consider another Drinfeld $A$-module $\phi'$ and an isogeny $f\colon\phi\to\phi'$ defined over~$K$. Recall that there exists an isogeny $g\colon \phi'\to\phi$ such that $g\circ f = \phi_a$ for some non-zero $a\in A$. From this it follows that $M' \defeq f(M)$ is a torsion free finitely generated $A$-submodule of $K$ for the action of $A$ through~$\phi'$. Thus $f$ induces $A_\ad$-linear maps from the modules in (\ref{TadM2}) to those associated to $\phi'$ and~$M'$. The existence of $g$ implies that these maps are inclusions of finite index. Together these maps yield a commutative diagram of $A_\ad$-modules with exact rows
\UseTheoremCounterForNextEquation
\begin{equation}\label{TadM2M'}
\vcenter{\xymatrix{ 
0 \ar[r] & T_\ad(\phi) \ar[r] \ar@{^{ (}->}[d] 
& T_\ad(\phi,M) \ar[r] \ar@{^{ (}->}[d] 
& M\otimes_AA_\ad \ar[r] \ar[d]^-{\wr}
& 0\\
0 \ar[r] & T_\ad(\phi') \ar[r] & T_\ad(\phi',M') \ar[r] 
& M'\otimes_AA_\ad \ar[r] & 0\rlap{.}\\}}
\end{equation}
% where the isomorphy on the right hand side results from the definition of~$M'$. 
By construction all these maps are equivariant under $\Gal(K^\sep/K)$; hence the images of Galois in each column are canonically isomorphic. If we denote the analogues of the groups $\Delta_{\ad,M}\subset\Gamma_{\ad,M}\onto\Gamma_\ad$ associated to $\phi'$ and $M'$ by $\Delta^{\phi'}_{\ad,M'}\subset\Gamma^{\phi'}_{\ad,M'}\onto\Gamma^{\phi'}_\ad$, this means that we have a natural commutative diagram
\UseTheoremCounterForNextEquation
\begin{equation}\label{DeltaAdMInclus2M'}
\vcenter{\xymatrix@C-13pt@R-5pt{ 
\Gamma_\ad && \ar@{->>}[ll] \Gamma_{\ad,M} & \ar@{}[l]|\supset
\Delta_{\ad,M} \ar@{^{ (}->}[rr] && \Hom_A(M,T_\ad(\phi)) \ar@{^{ (}->}[d] \\
\Gamma^{\phi'}_\ad \ar@{=}[u]_-{\wr} 
&& \ar@{->>}[ll] \Gamma^{\phi'}_{\ad,M'} \ar@{=}[u]_-{\wr} 
& \ar@{}[l]|\supset \Delta^{\phi'}_{\ad,M'} \ar@{=}[u]_-{\wr} \ar@{^{ (}->}[rr] 
&& \Hom_A(M',T_\ad(\phi')) \rlap{,}\\}}
\end{equation}
where the vertical arrow on the right hand side is an inclusion of finite index.

%%%%%%%%%%%%%%%%%%%%%%%%%%%%%%%%%%%%%%%%%%%%%%%

\medskip
Next let $B$ be any integrally closed infinite subring of~$A$. Then $A$ is a finitely generated projective $B$-module of some rank $s\ge1$. The restriction $\psi\defeq\phi|B$ is therefore a Drinfeld $B$-module of rank $rs$ over~$K$, and the given $A$-module $M$ of rank $d$ becomes a $B$-module of rank~$ds$. Moreover, since the characteristic of $\phi$ is by definition the kernel of the derivative map $a\mapsto d\phi_a$, the characteristic of $\psi$ is simply $\Fq_0 \defeq \Fp_0\cap B$. In analogy to (\ref{ATorDecomp2}) we have
\UseTheoremCounterForNextEquation
\begin{equation}\label{BTorDecomp2}
B_\ad\ \defeq\ \End_B(B_{(\Fq_0)}/B)
\ \cong\ \prod_{\Fq\not=\Fq_0} B_\Fq,
\end{equation}
where the product is extended over all maximal ideals $\Fq\not=\Fq_0$ of~$B$. Thus
\UseTheoremCounterForNextEquation
\begin{equation}\label{ABTorDecomp2}
A\otimes_BB_\ad\ \cong\ \prod_{\Fp\,\nmid\,\Fq_0} A_\Fp
\end{equation}
is in a natural way a factor ring of~$A_\ad$. More precisely, it is isomorphic to $A_\ad$ if the characteristic $\Fp_0$ and hence $\Fq_0$ is zero; otherwise it is obtained from $A_\ad$ by removing the finitely many factors $A_\Fp$ for all maximal ideals $\Fp\not=\Fp_0$ of $A$ above~$\Fq_0$. In particular we have a natural isomorphism $A\otimes_BB_\ad \cong A_\ad$ if and only if $\Fp_0$ is the unique prime ideal of $A$ above~$\Fq_0$.

\begin{Prop}\label{ChangingAProp}
The exact sequence (\ref{TadM2}) for $\psi$ and $M$ is naturally isomorphic to that obtained from the exact sequence (\ref{TadM2}) for $\phi$ and $M$ by tensoring with $A\otimes_BB_\ad$ over~$A_\ad$. In particular we have a commutative diagram with surjective vertical arrows
$$\xymatrix{ 
0 \ar[r] & T_\ad(\phi) \ar[r] \ar@{->>}[d] 
& T_\ad(\phi,M) \ar[r] \ar@{->>}[d] 
& M\otimes_AA_\ad \ar[r] \ar@{->>}[d] & 0\\
0 \ar[r] & T_\ad(\psi) \ar[r] & T_\ad(\psi,M) \ar[r] 
& M\otimes_BB_\ad \ar[r] & 0\rlap{.}\\}$$
If $\Fp_0$ is the only prime ideal of $A$ above~$\Fq_0$, the vertical arrows are isomorphisms.
\end{Prop}

\begin{Proof}
According to (\ref{DivHullPrimeDef}) the prime-to-$\Fp_0$ division hull of $M$ with respect to $\phi$ and the prime-to-$\Fq_0$ division hull of $M$ with respect to $\psi$ are
\begin{eqnarray*}
\Div_{K^\sep}^{(\Fp_0)}(M) &\defeq&
\bigl\{x\in K^\sep\;\bigm|\;\exists a\in A\setminus\Fp_0: \phi_a(x)\in M\bigr\},\\
\Div_{K^\sep}^{(\Fq_0)}(M) &\defeq&
\bigl\{x\in K^\sep\;\bigm|\;\exists b\in B\setminus\Fq_0: \psi_b(x)\in M\bigr\}.
\end{eqnarray*}
Here the latter is automatically contained in the former, because any $b\in B\setminus\Fq_0$ with $\psi_b(x)\in M$ is by definition an element $a := b\in A\setminus\Fp_0$ with $\phi_a(x)\in M$. Thus $\Div_{K^\sep}^{(\Fq_0)}(M)/M$ is the subgroup of all elements of $\Div_{K^\sep}^{(\Fp_0)}(M)/M$ that are annihilated by some element of $B\setminus\Fq_0$. In other words, it is the subgroup of all prime-to-$\Fq_0$ torsion with respect to~$B$, or again, it is obtained from $\Div_{K^\sep}^{(\Fp_0)}(M)/M$ by removing the $\Fp$-torsion for all maximal ideals $\Fp\not=\Fp_0$ of $A$ above~$\Fq_0$. In the same way $A\otimes_B(B_{(\Fq_0)}/B)$ is isomorphic to the submodule of $A_{(\Fp_0)}/A$ obtained by removing the $\Fp$-torsion for all $\Fp|\Fq_0$. The same process applied to the exact sequence (\ref{DivExDiag1}) therefore yields the analogue for $\psi$ and~$M$. By definition the exact sequence (\ref{TadM2}) for $\psi$ and $M$ is obtained from this by applying the functor 
$$X\ \mapsto\ \Hom_B(B_{(\Fq_0)}/B,X)
\ \cong\ \Hom_A\bigl(A\otimes_B(B_{(\Fq_0)}/B),X\bigr)$$
analogous to (\ref{TDef}). The total effect of this is simply to remove the $\Fp$-primary factors for all $\Fp|\Fq_0$ from the exact sequence (\ref{TadM2}) for $\phi$ and~$M$, from which everything follows.
\end{Proof}

By construction the diagram in Proposition \ref{ChangingAProp} is equivariant under $\Gal(K^\sep/K)$. It therefore induces a natural commutative diagram of Galois groups
\UseTheoremCounterForNextEquation
\begin{equation}\label{ChangingAGalois}
\vcenter{\xymatrix@C-13pt{ 
\Aut_{A_\ad}(T_\ad(\phi)) & \Gamma_\ad \ar@{}[l]|-\supset \ar@{->>}[d] 
&& \ar@{->>}[ll] \Gamma_{\ad,M} \ar@{->>}[d] & \ar@{}[l]|\supset
\Delta_{\ad,M} \ar@{^{ (}->}[rr] \ar[d] && \Hom_A(M,T_\ad(\phi)) \ar[d] \\
\Aut_{B_\ad}(T_\ad(\psi)) & \Gamma^\psi_\ad  \ar@{}[l]|-\supset
&& \ar@{->>}[ll] \Gamma^\psi_{\ad,M} & \ar@{}[l]|\supset
\Delta^\psi_{\ad,M} \ar@{^{ (}->}[rr] && \Hom_B(M,T_\ad(\psi)) \rlap{,}\\}}
\end{equation}
where the subgroups in the lower row are the analogues for $\psi$ and $M$ of those in the upper row. By construction the left two vertical arrows are surjective, and they are isomorphisms if $\Fp_0$ is the only prime ideal of $A$ above~$\Fq_0$. In that case the rightmost vertical arrow is injective and it follows that the map $\Delta_{\ad,M}\to\Delta^\psi_{\ad,M}$ is an isomorphism as well. In general one can only conclude that $\Delta^\psi_{\ad,M}$ contains the image of $\Delta_{\ad,M}$. In any case the diagram (\ref{ChangingAGalois}) gives a precise way of determining $\Gamma^\psi_{\ad,M}$ from $\Gamma_{\ad,M}$.

%%%%%%%%%%%%%%%%%%%%%%%%%%%%%%%%%%%%%%%%%%%%%%%%%%%%%%%%%%%%%%%%%%%%%%%%%%%%%%%%%%%%%%%%%%%%%

\section{Previous results on Galois groups}
\label{Glob1b}

% From now on we assume that $K$ is finitely generated over~$\BF_p$. 
% \medskip

In this section we recall some previous results on the Galois group~$\Gamma_\ad$. Its precise description up to commensurability depends on certain endomorphism rings. The endomorphism ring of a Drinfeld module of generic characteristic is always commutative, but in special characteristic it can be non-commutative. In the latter case it can grow on restricting $\phi$ to a subring $B$ of~$A$, and this effect can impose additional conditions on~$\Gamma_\ad$. The question of whether the endomorphism ring becomes stationary or grows indefinitely with $B$ depends on the following property:

\begin{Def}\label{IsotrivialDef}
We call a Drinfeld $A$-module of special characteristic over $K$ \emph{isotrivial} if over $K^\sep$ it is isomorphic to a Drinfeld $A$-module defined over a finite field.
\end{Def}

The next definition is slightly ad hoc, but it describes particular kinds of Drinfeld modules to which we can reduce ourselves in all cases, allowing a unified treatment of Kummer theory later on.

\begin{Def}\label{PrimitiveDef}
We call the triple $(A,K,\phi)$ \emph{primitive} if the following conditions hold:
\begin{enumerate}
\item[(a)] $R := \End_K(\phi)$ is equal to $\End_{K^\sep}(\phi)$.
\item[(b)] The center of~$R$ is~$A$.
\item[(c)] $R$ is a maximal $A$-order in $R\otimes_AF$.
\item[(d)] If $\phi$ is non-isotrivial of special characteristic, then for every integrally closed infinite subring $B\subset A$ we have $\End_{K^\sep}(\phi|B) = R$.
\item[(e)] If $\phi$ is isotrivial of special characteristic, then $A=\BF_p[a_0]$ with $\phi_{a_0} = \tau^{[k/\BF_p]}$, where $k$ denotes the finite field of constants of~$K$.
\end{enumerate}
\end{Def}

\begin{Prop}\label{PrimRed}
Let $A'$ denote the normalization of the center of $\End_{K^\sep}(\phi)$. 
\begin{enumerate}
\item[(a)] There exist a Drinfeld $A'$-module $\phi'\colon A'\to K^\sep[\tau]$ and an isogeny $f\colon\phi\to\phi'|A$ over~$K^\sep$ such that $A'$ is the center of $\End_{K^\sep}(\phi')$.
\item[(b)] The characteristic $\Fp'_0$ of any $\phi'$ as in (a) is a prime ideal of $A'$ above the characteristic $\Fp_0$ of~$\phi$.
\item[(c)] There exist a finite extension $K'\subset K^\sep$ of~$K$, a Drinfeld $A'$-module $\phi'\colon A'\to K'[\tau]$, an isogeny $f\colon\phi\to\phi'|A$ over~$K'$, and an integrally closed infinite subring $B\subset A'$ such that  $A'$ is the center of $\End_{K^\sep}(\phi')$ and $(B,K',\phi'|B)$ is primitive.
\item[(d)] The subring $B$ in (c) is unique unless $\phi$ is isotrivial of special characteristic, in which case it is never unique.
% not unique.
\item[(e)] For any data as in (c) the characteristic $\Fp'_0$ of $\phi'$ is the unique prime ideal of $A'$ above the characteristic $\Fq_0$ of~$\phi'|B$.
\end{enumerate}
\end{Prop}

\begin{Proof}
Applying \cite[Prop.$\,$4.3]{DevicPink} to $\phi$ and the center of $\End_{K^\sep}(\phi)$ yields a Drinfeld $A$-module $\tilde\phi\colon A\to K^\sep[\tau]$ and an isogeny $f\colon\phi\to\tilde\phi$ over $K^\sep$ such that $A'$ is mapped into $\End_{K^\sep}(\tilde\phi)$ under the isomorphism $\End_{K^\sep}(\phi)\otimes_AF \cong \End_{K^\sep}(\tilde\phi)\otimes_AF$ induced by~$f$. Then $A'\otimes_AF$ is the center of $\End_{K^\sep}(\tilde\phi)\otimes_AF$, and since $A'$ is integrally closed, it follows that $A'$ is the center of $\End_{K^\sep}(\tilde\phi)$. The tautological homomorphism $A' \into \End_{K^\sep}(\tilde\phi) \into K^\sep[\tau]$ thus constitutes a Drinfeld $A'$-module $\phi'$ with $\phi'|A \cong \tilde\phi$ and $\End_{K^\sep}(\phi')=\End_{K^\sep}(\tilde\phi)$. In particular the center of $\End_{K^\sep}(\phi')$ is equal to~$A'$, proving~(a).

For (b) recall that the characteristic of $\phi'$ is the kernel of the derivative map $a'\mapsto d\phi'_{a'}$. Calling it $\Fp_0'$, the characteristic of $\phi'|A$ is then $\Fp'_0\cap A$. As the characteristic of a Drinfeld module is invariant under isogenies, it follows that $\Fp'_0$ lies above~$\Fp_0$, proving~(b).

For the remainder of the proof we take any pair $\phi'$ and $f$ as in (a). We also choose a finite extension $K'\subset K^\sep$ of~$K$ such that $\phi'$ and $f$ are defined over $K'$ and that $\End_{K'}(\phi') = \End_{K^\sep}(\phi')$.

Suppose first that $\phi'$ has generic characteristic. Then $\End_{K'}(\phi')$ is commutative and hence equal to~$A'$, and so the triple $(A',K',\phi')$ is already primitive. This proves (c) with $B=A'$. Also, for any integrally closed infinite subring $B\subset A'$ the ring $\End_{K^\sep}(\phi'|B)$ is commutative and hence again equal to $\End_{K^\sep}(\phi')=A'$. Thus $(B,K',\phi'|B)$ being primitive requires that $B=A'$, proving~(d). Since $B=A'$, the assertion of (e) is then trivially true.

Suppose next that $\phi'$ is non-isotrivial of special characteristic. Then by \cite[Thm.$\,$6.2]{PinkDrinSpec2} there exists a unique integrally closed infinite subring $B\subset A'$ such that $B$ is the center of $\End_{K^\sep}(\phi'|B)$ and that $\End_{K^\sep}(\phi'|B') \subset \End_{K^\sep}(\phi'|B)$ for every integrally closed infinite subring $B'\subset A'$. 
For use below we note that both properties are invariant under isogenies of~$\phi'$, because isogenies induce isomorphisms on the rings $\End_{K^\sep}(\phi'|B')\otimes_{B'}\Quot(B')$. By uniqueness it follows that~$B$, too, is invariant under isogenies of~$\phi'$. 

After replacing $K'$ by a finite extension we may assume that $\End_{K'}(\phi'|B) = \End_{K^\sep}(\phi'|B)$. Then the triple $(B,K',\phi'|B)$ satisfies the conditions in Definition \ref{PrimitiveDef} except that $\End_{K'}(\phi'|B)$ may not be a maximal order in $\End_{K'}(\psi')\otimes_B\Quot(B)$. But applying \cite[Prop.$\,$4.3]{DevicPink} to $\phi'|B$ and $\End_{K'}(\phi'|B)$ yields a Drinfeld $B$-module $\psi'\colon B\to K'[\tau]$ and an isogeny $g\colon\phi'|B\to\psi'$ over $K'$ such that $\End_{K'}(\psi')$ is a maximal order in $\End_{K'}(\psi')\otimes_B\Quot(B)$ which contains $\End_{K'}(\phi'|B)$. By the preceding remarks we now find that $(B,K',\psi')$ is primitive. 
Moreover, the composite homomorphism $A' \into \End_{K'}(\phi') \into \End_{K'}(\phi'|B) \into \End_{K'}(\psi') \into K'[\tau]$ constitutes a Drinfeld $A'$-module $\phi''$ with $\phi''|B \cong \psi'$. After replacing $(\phi',f)$ by $(\phi''\!,g\:{\circ}f)$ the data then satisfies all the requirements of (c). 
Assertion (d) follows from the above stated uniqueness of~$B$, and (e) follows from \cite[Prop.$\,$3.5]{PinkDrinSpec2}. 

It remains to consider the case where $\phi'$ is isotrivial of special characteristic. In this case we may assume that $\phi'$ is defined over the constant field $k'$ of~$K'$. 
% factors through $k'[\tau]\subset K'[\tau]$ where $k'$ denotes the constant field~$K'$. 
Any endomorphism of $\phi'$ over $K^\sep$ is then defined over a finite extension of~$k'$, but by assumption also over~$K'$; hence it is defined over~$k'$. In other words we have $\End_{K^\sep}(\phi') \subset k'[\tau]$. Since $\tau^{[k'/\BF_p]}$ lies in the center of the $k'[\tau]$, it thus corresponds to an element of the center of $\End_{K^\sep}(\phi')$. As this center is equal to $A'$ by assumption, there is therefore an element $a_0\in A'$ with $\phi'_{a_0} = \tau^{[k'/\BF_p]}$. Set $B \defeq \BF_p[a_0] \subset A'$ which, being isomorphic to a polynomial ring, is an integrally closed infinite subring of~$A'$. Then $\End_{K^\sep}(\phi'|B)$ is the commutant of $\tau^{[k'/\BF_p]}$ in $K^\sep[\tau]$ and hence just $k'[\tau]$. By a standard construction this is a maximal $B$-order in a (cyclic) central division algebra over $\Quot(B)$; hence $(B,K',\phi'|B)$ is primitive, proving (c). 
For (d) observe that replacing $K'$ and $k'$ by finite extensions amounts to replacing $a_0$ by an arbitrary positive power~$a_0^i$. Thus the ring $B$ is really not unique in this case, proving (d). Finally, assertion (e) follows from \cite[Prop.$\,$6.4$\,$(a)]{DevicPink}.
This finishes the proof of Proposition \ref{PrimRed}.
\end{Proof}

%%%%%%%%%%%%%%%%%%%%%%%%%%%%%%%%%%%%%%%%%%%%%%%

Assume now that $(A,K,\phi)$ is primitive and that $\phi$ has rank~$r$. Then $R\otimes_AF$ is a central division algebra of dimension $m^2$ over $F$ for some factorization $r=mn$. Thus for all primes $\Fp\not=\Fp_0$ of~$A$, the ring $R_\Fp := R\otimes_AA_\Fp$ is an $A_\Fp$-order in the central simple algebra $R\otimes_AF_\Fp$ of dimension $m^2$ over $F_\Fp$ and is isomorphic to the matrix ring $\Mat_{m\times m}(A_\Fp)$ for almost all~$\Fp$. Let $D_\Fp$ denote the commutant of $R_\Fp$ in $\End_{A_\Fp}(T_\Fp(\phi))$. This is an $A_\Fp$-order in a central simple algebra of dimension $n^2$ over~$F_\Fp$ and is isomorphic to the matrix ring $\Mat_{n\times n}(A_\Fp)$ for almost all~$\Fp$. Let $D^1_\Fp$ denote the multiplicative group of elements of $D_\Fp$ of reduced norm~$1$. Set 
$$D_\ad\ :=\ \prod_{\Fp\not=\Fp_0} D_\Fp\ \subset\ \End_{A_\ad}(T_\ad(\phi))$$
and 
$$D^1_\ad\ :=\ \prod_{\Fp\not=\Fp_0} D^1_\Fp \ \subset\ D_\ad^\times
\ \subset\ \Aut_{A_\ad}(T_\ad(\phi)).$$
If $\phi$ has generic characteristic, we have $m=1$ and therefore $D_\Fp = \End_{A_\Fp}(T_\Fp(\phi)) \cong \Mat_{r\times r}(A_\Fp)$ for all~$\Fp$.

If $\phi$ is non-isotrivial of special characteristic~$\Fp_0$, let $a_0$ be any element of $A$ that generates a positive power of~$\Fp_0$. If $\phi$ is isotrivial, the element $a_0$ from Definition \ref{PrimitiveDef} (d) already has the same property. In both cases we view $a_0$ as a scalar element of $D_\ad^\times$ via the diagonal embedding $A\subset A_\ad\subset D_\ad$, and let $\smash{\overline{\langle a_0\rangle}}$ denote the pro-cyclic subgroup that is topologically generated by it. 

\medskip
In general the group $\Gamma_\ad$ was described up to commensurability in Pink-R\"utsche \cite{PinkRuetsche2} and Devic-Pink \cite{DevicPink}. In the primitive case \cite[Thm.$\;$0.1]{PinkRuetsche2} and \cite[Thm.$\;$1.1, Prop.$\,$6.3]{DevicPink} imply:

\begin{Thm}\label{GammaAdelicOpen}
Assume that $(A,K,\phi)$ is primitive.
\begin{enumerate}
\item[(a)] If $\phi$ has generic characteristic, then $\Gamma_\ad$ is open in $D_\ad^\times$.
\item[(b)] If $\phi$ is non-isotrivial of special characteristic, then $n\ge2$ and $\Gamma_\ad$ is commensurable with $\overline{\langle a_0\rangle} \cdot D^1_\ad$.
\item[(c)] If $\phi$ is isotrivial of special characteristic, then $n=1$ and $\Gamma_\ad = \overline{\langle a_0\rangle}$ with $a_0$ from \ref{PrimitiveDef}~(c).
\end{enumerate}
\end{Thm}

\begin{Cor}\label{GammaIdeals}
Assume that $(A,K,\phi)$ is primitive.
\begin{enumerate}
\item[(a)] Let $\Theta_\ad$ denote the closure of the $\BF_p$-subalgebra of $D_\ad$ generated by~$\Gamma_\ad$. Then there exists a non-zero ideal $\Fa$ of $A$ with $\Fa\not\subset\Fp_0$ such that $\Fa D_\ad\subset\Theta_\ad$.
\item[(b)] There exist a scalar element $\gamma\in\Gamma_\ad$ and a non-zero ideal $\Fb$ of $A$ with $\Fb\not\subset\Fp_0$ such that $\gamma\equiv 1$ modulo $\Fb A_\ad$ but not modulo $\Fp\Fb A_\ad$ for any prime $\Fp\not=\Fp_0$ of~$A$.
\end{enumerate}
\end{Cor}

\begin{Proof}
Any open subgroup of $D_\ad^\times$, and for $n\ge2$ any open subgroup of~$D_\ad^1$, generates an open subring of~$D_\ad$. Thus the assertion (a) follows from Theorem \ref{GammaAdelicOpen} unless $\phi$ is isotrivial of special characteristic. But in that case we have $\Theta_\ad = \smash{\overline{\BF_p[a_0]}} = A_\ad = D_\ad$ and (a) follows as well.

In generic characteristic the assertion (b) follows directly from the openness of~$\Gamma_\ad$. In special characteristic some positive power $a_0^i$ lies in $\Gamma_\ad$, and so (b) holds with $\gamma=a_0^i$ and the ideal $\Fb = (a_0^i-1)$.
\end{Proof}

%%%%%%%%%%%%%%%%%%%%%%%%%%%%%%%%%%%%%%%%%%%%%%%%%%%%%%%%%%%%%%%%%%%%%%%%%%%%%%%%%%%%%%%%%%%%%

\section{The primitive case}
\label{PrimCase}

% Keeping the notations and assumptions of the preceding sections, we 
Now we prove the following result, of which Theorem \ref{MainIntro} (c) is a special case:

\begin{Thm}\label{MainPrimitive}
Assume that $(A,K,\phi)$ is primitive. Set $R:= \End_K(\phi)$ and let $M$ be a finitely generated torsion free $R$-submodule of~$K$. Then the inclusions $\Delta_{\ad,M} \subset \Hom_R(M,T_\ad(\phi))$ and $\Gamma_{\ad,M} \subset \Gamma_\ad \ltimes \Hom_R(M,T_\ad(\phi))$ are both open.
\end{Thm}

So assume that $(A,K,\phi)$ is primitive. Let the subring $\Theta_\ad \subset D_\ad$, the element $\gamma\in\Gamma_\ad$, and the ideals $\Fa$, $\Fb\subset A$ be as in Corollary \ref{GammaIdeals}. Since $M \subset \Div_K^{(\Fp_0)}(M)$ has finite index by Theorem \ref{lettThm}, we can also choose a non-zero ideal $\Fc$ of $A$ with $\Fc\not\subset\Fp_0$ such that $\Fc\cdot \Div_K^{(\Fp_0)}(M) \subset M$. With this data we prove the following more precise version of Theorem \ref{MainPrimitive}:

\begin{Thm}\label{MainPrimSpecific}
In the above situation we have $\Fa\Fb\Fc \cdot\Hom_R(M,T_\ad(\phi)) \subset \Delta_{\ad,M}$.
\end{Thm}

%%%%%%%%%%%%%%%%%%%%%%%%%%%%%%%%%%%%%%%%%%%%%%%

In the rest of this section we abbreviate $T_\ad := T_\ad(\phi)$ and $M^*_\ad \defeq \Hom_R(M,T_\ad)$. Recall that the embedding $\Delta_{\ad,M} \subset M^*_\ad$ is adjoint to the pairing $\langle\ ,\ \rangle$ from (\ref{PairingDef}). The arithmetic part of the proof is a calculation in $\Div_{K^\sep}^{(\Fp_0)}(M)$ with the following result:

\begin{Lem}\label{MainPrim1}
For any prime $\Fp\not=\Fp_0$ of $A$ and any element $m\in M$ satisfying $\langle \Delta_{\ad,M}, m\rangle \subset \Fp\Fb\Fc T_\ad$ we have $m\in\Fp M$.
\end{Lem}

\begin{Proof}
The assumption on $\gamma$, viewed as an element of~$A_\ad$, means that $\gamma-1\in \Fb A_\ad \setminus \Fp\Fb A_\ad$. By the Chinese remainder theorem we can find an element $b\in A\setminus\Fp_0$ satisfying $b\equiv\gamma-1$ modulo $\Fp\Fb\Fc A_\ad$. Then by construction we have $b\in\Fb\setminus\Fp\Fb$, and $\gamma$ acts on all $\Fp\Fb\Fc$-torsion points of $\phi$ through the action of $1+b\in A$.
By the Chinese remainder theorem we can also find elements $a$, $c\in A\setminus\Fp_0$ with $a\in\Fp\setminus\Fp^2$ and $c\in\Fc\setminus\Fp\Fc$. Then the product $abc$ lies in $\Fp\Fb\Fc\setminus(\Fp^2\Fb\Fc\cup\Fp_0)$. In particular the order of $abc$ at $\Fp$ is equal to that of $\Fp\Fb\Fc$, and so we can also choose an element $d\in A\setminus\Fp$ such that $d\Fp\Fb\Fc \subset (abc)$.

For better readability we abbreviate the action of any element $e\in A$ on an element $x\in K^\sep$ by $ex \defeq \phi_e(x)$. 
% This should not be too confusing, because $\phi$ is fixed.
Since $abc\in A\setminus\Fp_0$, we can select an element $\Tildem \in \Div_{K^\sep}^{(\Fp_0)}(M)$ with $abc\Tildem = m$. Then $\tildem \defeq d\Tildem$ is an element of $\Div_{K^\sep}^{(\Fp_0)}(M)$ which satisfies $abc\tildem = dm$. By construction $\tildem$ lies in~$K^\sep$, but we shall see that it actually lies in a specific subfield. 

Namely, choosing a compatible system of division points of $\tildem$ we can find an $A$-linear map $\tilde t\colon A_{(\Fp_0)} \to \Div_{K^\sep}^{(\Fp_0)}(M)$ satisfying $\tilde t\bigl(\frac{1}{abc}\bigr) = \tildem$. Then $\tilde t(1) = abc\tilde m = dm$ lies in~$M$; hence $\tilde t$ induces an $A$-linear map $t\colon A_{(\Fp_0)}/A \to \Div_{K^\sep}^{(\Fp_0)}(M)/M$. By the construction (\ref{TadMDef}) this map is an element of $T_\ad(\phi,M)$ whose image in ${M\otimes_AA_\ad}$ is $dm\otimes1$. 
For any $\delta\in\Delta_{\ad,M}$ the definition (\ref{PairingDef}) of the pairing now says that $\langle\delta,dm\rangle = \delta(t)-t$.
But the assumption $\langle \Delta_{\ad,M}, m\rangle \subset \Fp\Fb\Fc T_\ad$ implies that $\langle \Delta_{\ad,M}, dm \rangle \subset d\Fp\Fb\Fc T_\ad \allowbreak \subset \nobreak abc T_\ad$ and therefore $\delta(t)-t \in abc T_\ad$.
Thus $\delta(t)-t$ is the multiple by $abc$ of an $A$-linear map $A_{(\Fp_0)}/A \to \Div_{K^\sep}^{(\Fp_0)}(\{0\})$; hence it is zero on the residue class of~$\frac{1}{abc}$.
By the construction of $t$ this means that $\delta(\tildem)-\tildem=0$.
Varying $\delta$ we conclude that $\tildem$ is fixed by $\Delta_{\ad,M}$; in other words, it lies in the subfield $K_\ad\subset K^\sep$ with $\Gal(K_\ad/K)=\Gamma_\ad$.

Now consider any element $\sigma\in \Gal(K^\sep/K)$. The fact that $m$ lies in $K$ implies that
$$abc(\sigma-1)(\Tildem) = (\sigma-1)(abc\Tildem) = (\sigma-1)(m) = 0.$$
Thus $(\sigma-1)(\Tildem)$ is annihilated by $abc$ and hence by the ideal $d\Fp\Fb\Fc \subset (abc)$. The element $(\sigma-1)(\tildem) = d(\sigma-1)(\Tildem)$ is therefore annihilated by the ideal $\Fp\Fb\Fc$. 
% Note: It's for this step that we need to calculate with the element $\Tildem$!
Since $\gamma$ acts on all $\Fp\Fb\Fc$-torsion points through the action of $1+b\in A$, it follows that 
$$(1+b-\gamma)((\sigma-1)(\tildem)) = 0.$$
On the other hand we have $\tildem\in K_\ad$, and since $\gamma$ lies in the center of $\Gamma_\ad$, its action on $K_\ad$ commutes with the action of $\sigma$ on~$K_\ad$. Thus the last equation is equivalent to
$$(\sigma-1)((1+b-\gamma)(\tildem)) = 0.$$
As $\sigma\in \Gal(K^\sep/K)$ was arbitrary, it follows that $(1+b-\gamma)(\tildem) \in K$.

Since $\tildem$ lies in $\Div_{K^\sep}^{(\Fp_0)}(M)$, we can now deduce that $(1+b-\gamma)(\tildem) \in \Div_K^{(\Fp_0)}(M)$. By the choice of $\Fc$ and $c$ it follows that $c(1+b-\gamma)(\tildem) \in M$. The fact that $dm=abc\tildem$ thus implies that
$$(1+b-\gamma)(dm) =
(1+b-\gamma)(abc\tildem) =
abc(1+b-\gamma)(\tildem) \in abM.$$
But $dm$ is an element of~$K$ and therefore satisfies $(1-\gamma)(dm)=0$. Thus the last relation shows that actually $bdm\in abM$ and so $dm\in aM$. Since $a\in\Fp$ and $d\not\in\Fp$, this implies that $m\in\Fp M$, as desired.
\end{Proof}

%%%%%%%%%%%%%%%%%%%%%%%%%%%%%%%%%%%%%%%%%%%%%%%

\medskip
The rest of the proof of Theorem \ref{MainPrimSpecific} involves rings and modules. It is easier to understand if $R=A$, in which case $D_\Fp\cong \Mat_{r\times r}(A_\Fp)$ for all~$\Fp$, so the readers may want to restrict themselves to that case on first reading. For general facts on maximal orders in semisimple algebras see Reiner \cite{ReinerMaximalOrders}.

\medskip
Using \cite[Cor.$\,$11.6]{ReinerMaximalOrders}, the assumptions \ref{PrimitiveDef} imply that for any prime $\Fp\not=\Fp_0$ of $A$ the ring $R_\Fp := R\otimes_AA_\Fp$ is a maximal order in a finite dimensional central simple algebra over~$F_\Fp$. 
By \cite[Thm.$\,$17.3]{ReinerMaximalOrders} we can therefore identify it with the matrix ring $\Mat_{n_\Fp\times n_\Fp}(S_\Fp)$, where $S_\Fp$ is the maximal order in a finite dimensional central division algebra over~$F_\Fp$. 
Here $n_\Fp\ge1$ may vary with~$\Fp$. 
Let $L_\Fp := S_\Fp^{\oplus n_\Fp}$ denote the tautological left $R_\Fp$-module. Since $T_\Fp := T_\Fp(\phi)$ is a non-trivial finitely generated torsion free left $R_\Fp$-module, it is isomorphic to $L_\Fp^{\oplus m_\Fp}$ for some $m_\Fp\ge1$ by \cite[Thm.$\,$18.10]{ReinerMaximalOrders}.
Thus $D_\Fp \defeq \End_{R_\Fp}(T_\Fp)$ is isomorphic to the matrix ring $\Mat_{m_\Fp\times m_\Fp}(S^\opp_\Fp)$ over the opposite algebra $S^\opp_\Fp$. Let $N_\Fp := (S_\Fp^\opp)^{\oplus m_\Fp}$ denote the tautological left $D_\Fp$-module; then as a $D_\Fp$-module $T_\Fp$ is isomorphic to $N_\Fp^{\oplus n_\Fp}$. Moreover, by biduality, using Morita equivalence \cite[Thm.$\,$16.14]{ReinerMaximalOrders} or direct computation, we have
\UseTheoremCounterForNextEquation
\begin{equation}\label{Biduality1}
R_\Fp \cong \End_{D_\Fp}(T_\Fp).
\end{equation}

Next, since $R$ is a maximal order in a division algebra over the Dedekind ring~$A$, and $M$ is a finitely generated torsion free $R$-module, $M$ is a projective $R$-module by \cite[Cor.$\,$21.5]{ReinerMaximalOrders}, say of rank $\ell\ge0$.
For each $\Fp\not=\Fp_0$ we therefore have $M_\Fp \defeq M\otimes_AA_\Fp \cong R_\Fp^{\oplus\ell}$ as an $R_\Fp$-module. Consequently $M_\Fp^* \defeq \Hom_R(M,T_\Fp) \cong T_\Fp^{\oplus\ell}$ as a $D_\Fp$-module via the action of $D_\Fp$ on~$T_\Fp$. Using the biduality (\ref{Biduality1}) we obtain a natural isomorphism
\UseTheoremCounterForNextEquation
\begin{equation}\label{Biduality2}
M_\Fp \cong \Hom_{D_\Fp}(M_\Fp^*,T_\Fp).
\end{equation}
Taking the product over all $\Fp\not=\Fp_0$ yields adelic versions of all this with $T_\ad = \prod T_\Fp$ and $R_\ad = \prod R_\Fp$ and $D_\ad = \prod D_\Fp$ and $M_\ad^* = \prod M_\Fp^*$.

\medskip
Recall that $\Delta_{\ad,M}$ is a closed additive subgroup of $M_\ad^* = \prod M_\Fp^*$. Let $\Delta_\Fp$ denote its image under the projection to~$M_\Fp^*$. 

\begin{Lem}\label{MainPrim2}
For any  $\Fp\not=\Fp_0$ and any $D_\Fp$-linear map $f\colon M_\Fp^* \to N_\Fp$ satisfying $f(\Delta_\Fp) \subset \Fp\Fb\Fc N_\Fp$, we have $f(M_\Fp^*) \subset \Fp N_\Fp$.
\end{Lem}

\begin{Proof}
Since $N_\Fp$ is a $D_\Fp$-module isomorphic to a direct summand of~$T_\Fp$, it is equivalent to show that for every $D_\Fp$-linear map $g\colon M_\Fp^* \to T_\Fp$ with $g(\Delta_\Fp) \subset \Fp\Fb\Fc T_\Fp$ we have $g(M_\Fp^*) \subset \Fp T_\Fp$. Let $\langle\ ,\ \rangle\colon \allowbreak M^*_\Fp\times M_\Fp \to T_\Fp$ denote the natural $A_\Fp$-bilinear map.
Then the biduality (\ref{Biduality2}) says that $g=\langle\underline{\ \ },m_\Fp\rangle$ for 
% any such map $g$ is represented by 
an element $m_\Fp\in M_\Fp$. Write $\Fp\Fb\Fc = \Fp^i\Fd$ for an integer $i\ge1$ and an ideal $\Fd$ of $A$ that is prime to~$\Fp$. Choose any element $m\in M$ which is congruent to $m_\Fp$ modulo $\Fp^i M_\Fp$ and congruent to $0$ modulo $\Fd M$. 
% Then the map $g$ is congruent to the map $\langle\ ,m\rangle$ modulo
Then the assumption 
$\langle \Delta_\Fp,m_\Fp\rangle = g(\Delta_\Fp) \subset \Fp\Fb\Fc T_\Fp$ implies that $\langle \Delta_{\ad,M}, m\rangle \subset \Fp\Fb\Fc T_\ad$. By Lemma \ref{MainPrim1} it follows that $m\in\Fp M$. Consequently $m_\Fp\in \Fp M_\Fp$ and therefore $g(M_\Fp^*) = \langle \Delta_\Fp, m_\Fp\rangle \subset \Fp T_\Fp$, as desired.
\end{Proof}

Now observe that $\Delta_{\ad,M} \subset M_\ad^*$ is a closed additive subgroup that is invariant under the action of~$\Gamma_\ad$. It is therefore a submodule with respect to the subring $\Theta_\ad \defeq \smash{\overline{\BF_p[\Gamma_\ad]}}$ of $D_\ad$ from Corollary \ref{GammaIdeals} (a). By \ref{GammaIdeals} (a) we therefore have
\UseTheoremCounterForNextEquation
\begin{equation}\label{DeltaPrime}
\Delta'_\ad \defeq \Fa D_\ad\Delta_{\ad,M} \subset\Delta_{\ad,M}.
\end{equation}
By construction $\Delta'_\ad$ is a submodule over $D_\ad = \prod D_\Fp$ and therefore itself a product $\Delta'_\ad = \prod\Delta'_\Fp$ for $D_\Fp$-submodules $\Delta'_\Fp \subset M^*_\Fp$.

\begin{Lem}\label{MainPrim3}
For any  $\Fp\not=\Fp_0$ and any $D_\Fp$-linear map $f\colon M_\Fp^* \to N_\Fp$ satisfying $f(\Delta'_\Fp) \subset \Fp\Fa\Fb\Fc N_\Fp$, we have $f(M_\Fp^*) \subset \Fp N_\Fp$.
\end{Lem}

\begin{Proof}
The definition of $\Delta'_\Fp$ implies that $\Fa D_\Fp\Delta_\Fp \subset\Delta'_\Fp$. Thus by assumption we have $\Fa f(\Delta_\Fp) \subset f(\Fa D_\Fp\Delta_\Fp) \subset f(\Delta'_\Fp) \subset \Fp\Fa\Fb\Fc N_\Fp$ and therefore $f(\Delta_\Fp) \subset \Fp\Fb\Fc N_\Fp$. By Lemma \ref{MainPrim2} this implies that $f(M_\Fp^*) \subset \Fp N_\Fp$.
\end{Proof}

\begin{Lem}\label{MainPrim4}
For any $\Fp\not=\Fp_0$ we have $\Fa\Fb\Fc M_\Fp^* \subset \Delta'_\Fp$.
\end{Lem}

\begin{Proof}
Let $\Fm_\Fp$ denote the maximal ideal of $S_\Fp^\opp$. Then by \cite[Thm.$\,$13.2]{ReinerMaximalOrders} we have $\Fp S_\Fp^\opp = \Fm_\Fp^e$ for some integer $e\ge1$. 
The general theory says the following about the structure of the module $M_\Fp^*/\Delta'_\Fp$ over the maximal order~$D_\Fp$. On the one hand, by Knebusch \cite[Satz$\,$7]{Knebusch1967} the torsion submodule of $M_\Fp^*/\Delta'_\Fp$ is a finite direct sum of indecomposable modules isomorphic to $\smash{N_\Fp/\Fm_\Fp^{j_\nu} N_\Fp}$ for certain integers $j_\nu\ge1$. On the other hand, the factor module of $M_\Fp^*/\Delta'_\Fp$ by its torsion submodule is projective by \cite[Cor.$\,$21.5]{ReinerMaximalOrders} and hence isomorphic to a direct sum of copies of~$N_\Fp$. That the factor module is projective also implies that $M_\Fp^*/\Delta'_\Fp$ is isomorphic to the direct sum of its torsion submodule with the factor module. Together it follows that $M_\Fp^*/\Delta'_\Fp$ is a finite direct sum of modules isomorphic to $N_\Fp$ or to $\smash{N_\Fp/\Fm_\Fp^{j_\nu} N_\Fp}$ for certain integers $j_\nu\ge1$. 

To use this fact, let $\Fp^i$ denote the highest power of $\Fp$ dividing $\Fa\Fb\Fc$.
If no summand isomorphic to $N_\Fp$ occurs in $M_\Fp^*/\Delta'_\Fp$ and all exponents $j_\nu$ are $\le ei$, then $M_\Fp^*/\Delta'_\Fp$ is annihilated by $\Fp^i S_\Fp^\opp = \Fm_\Fp^{ei}$. In this case it follows that $\Fa\Fb\Fc M_\Fp^* = \Fp^iM_\Fp^* \subset \Delta'_\Fp$, as desired. 

Otherwise there exists a surjective $D_\Fp$-linear map $M_\Fp^*/\Delta'_\Fp \onto N_\Fp/\Fm_\Fp^{ei+1} N_\Fp$. Composed with the isomorphism 
$$N_\Fp/\Fm_\Fp^{ei+1} N_\Fp \ \cong\  
\Fm_\Fp^{e-1}N_\Fp/\Fp^{i+1} N_\Fp \ =\ 
\Fm_\Fp^{e-1}N_\Fp/\Fp\Fa\Fb\Fc N_\Fp,$$
this yields a $D_\Fp$-linear map $M_\Fp^*/\Delta'_\Fp \to N_\Fp/\Fp\Fa\Fb\Fc N_\Fp$ whose image is not contained in $\Fp N_\Fp/\Fp\Fa\Fb\Fc N_\Fp = \Fm_\Fp^e N_\Fp/\Fp\Fa\Fb\Fc N_\Fp$. As $M_\Fp^*$ is a projective $D_\Fp$-module, the latter map can be lifted to a $D_\Fp$-linear map $f\colon M_\Fp^* \to N_\Fp$. By construction this map then satisfies $f(\Delta'_\Fp) \subset \Fp\Fa\Fb\Fc N_\Fp$ and $f(M_\Fp^*) \not\subset \Fp N_\Fp$. But that contradicts Lemma \ref{MainPrim3}; hence this case is not possible and the lemma is proved.
\end{Proof}

Taking the product over all $\Fp$, Lemma \ref{MainPrim4} and the inclusion (\ref{DeltaPrime}) imply that $\Fa\Fb\Fc M_\ad^* \subset \Delta'_\ad \subset\Delta_{\ad,M}$. This finishes the proof of Theorem \ref{MainPrimSpecific}. In particular it proves the first assertion of Theorem \ref{MainPrimitive}, from which the second assertion directly follows.

%%%%%%%%%%%%%%%%%%%%%%%%%%%%%%%%%%%%%%%%%%%%%%%%%%%%%%%%%%%%%%%%%%%%%%%%%%%%%%%%%%%%%%%%%%%%%

\section{The general case}
\label{GenCase}

First we note the following general fact on homomorphisms of modules:

\begin{Prop}\label{QuasiHom}
Let $S$ be a unitary ring, not necessarily commutative, and let $M$ and $N$ be left $S$-modules. Let $X$ be a subset of $M$ and $SX$ the $S$-submodule generated by it. Let $\Hom_{(S)}(X,N)$ denote the set of maps $\ell\colon X\to N$ such that for any finite collection of $s_i\in S$ and $x_i\in X$ with $\sum_i s_ix_i=0$ in $M$ we have $\sum_i s_i\ell(x_i)=0$ in~$N$.
\begin{enumerate}
\item[(a)] The restriction of maps induces a bijection $\Hom_{S}(SX,N) \stackrel{\sim}{\to}\Hom_{(S)}(X,N)$.
\item[(b)] If $R$ is a unitary subring of $S$ such that $X$ is an $R$-submodule and the natural map 
% from the non-commutative tensor product 
$S\otimes_RX\to M$, ${\sum_i s_i\otimes x_i} \mapsto \sum_i s_ix_i$ is injective, then $\Hom_{(S)}(X,N) = \Hom_{R}(X,N)$.
\item[(c)] If $X$ is an $S$-submodule of $M$, then $\Hom_{(S)}(X,N) = \Hom_{S}(X,N)$.
\end{enumerate}
\end{Prop}

\begin{Proof}
Let $F := \bigoplus_{x\in X} S\cdot[x]$ be the free left $S$-module over the set~$X$ and consider the natural $S$-linear map $F\to M$, $\sum_i s_i[x_i] \mapsto \sum_i s_ix_i$. Since $S$ is unitary, the image of this map is~$SX$. Let $T$ denote its kernel. Then giving an $S$-linear map $SX\to N$ is equivalent to giving an $S$-linear map $F\to N$ which vanishes on~$T$. Using the universal property of $F$ we find that the latter is equivalent to giving an element of $\Hom_{(S)}(X,N)$. The total correspondence is given by restriction of maps, proving~(a).

In (b) we have $S\otimes_RX\stackrel{\sim}{\to}SX$; hence the adjunction between tensor product and $\Hom$ yields bijections $\Hom_{S}(SX,N) \stackrel{\sim}{\to} \Hom_{S}(S\otimes_RX,N) \stackrel{\sim}{\to} \Hom_{R}(X,N)$. Their composite is again just restriction of maps; so by (a) the restriction map $\Hom_{(S)}(X,N) \to \Hom_{R}(X,N)$ is also bijective, proving~(b).

Finally, (c) is a special case of (a) or (b), according to taste.
\end{Proof}

%%%%%%%%%%%%%%%%%%%%%%%%%%%%%%%%%%%%%%%%%%%%%%%

\medskip
Now we return to the situation of Section~\ref{Glob1b}. We choose data $(A',K',\phi',f,B)$ as in Proposition \ref{PrimRed} (c), that is: We let $A'$ denote the normalization of the center of $\End_{K^\sep}(\phi)$, take a finite extension $K'\subset K^\sep$ of~$K$, a Drinfeld $A'$-module $\phi'\colon A'\to K'[\tau]$, an isogeny $f\colon\phi\to\phi'|A$ over~$K'$, and an integrally closed infinite subring $B\subset A'$ such that  $A'$ is the center of $\End_{K^\sep}(\phi')$ and $(B,K',\phi'|B)$ is primitive. By Proposition \ref{PrimRed} (e) the characteristic $\Fp'_0$ of $\phi'$ is then the only prime ideal of $A'$ above the characteristic $\Fq_0$ of~$\phi'|B$. We will apply the reduction steps from Section \ref{RedStep} to the isogeny $f$ and to each of the inclusions $A\subset A'\supset B$.

Specifically, let us set $\psi' \defeq \phi'|B$ and $S' \defeq \End_{K^\sep}(\psi')$. Then $M' \defeq A'f(M)$ is a finitely generated $A'$-submodule of $K'$ for the action of $A'$ through~$\phi'$, and so $N' \defeq S'M'$ is a finitely generated $B$-submodule for the action of $B$ through~$\psi'$. The modules $M'$ and $N'$ may have torsion, but since they are finitely generated, their torsion is annihilated by some non-zero element $a'\in A'$. Replacing the isogeny $f$ by $\phi'_{a'}\circ f$ replaces $M'$ by $a'M'$ and $N'$ by~$a'N'$; hence we may without loss of generality assume that $M'$ and $N'$ are torsion free.

Let $\Delta^{\phi'}_{\ad,M'}\subset\Gamma^{\phi'}_{\ad,M'}\onto\Gamma^{\phi'}_\ad$ denote the Galois groups as in (\ref{DelGamGam}) associated to $(K',\phi',M')$ in place of $(K,\phi,M)$, and similarly for $(K',\phi',N')$, respectively for $(K',\psi',N')$, and so on.
Then Proposition \ref{ChangingAProp} for the inclusion $A'\supset B$ yields a natural commutative diagram with exact rows
\UseTheoremCounterForNextEquation
\begin{equation}\label{ChangingAProp1}
\vcenter{\xymatrix@R-5pt{ 
0 \ar[r] & T_\ad(\phi') \ar[r] \ar@{=}[d]^-{\wr}
& T_\ad(\phi',N') \ar[r] \ar@{=}[d]^-{\wr}
& N'\otimes_{A'}A'_\ad \ar[r] \ar@{=}[d]^-{\wr} & 0\\
0 \ar[r] & T_\ad(\psi') \ar[r] & T_\ad(\psi',N') \ar[r] 
& N'\otimes_BB_\ad \ar[r] & 0\rlap{.}\\}}
\end{equation}
The action of $S'$ on the lower row thus yields a natural action on the upper row. Recall from (\ref{GammaAdMInclus2}) that any $S'$-equivariant splitting induces an embedding
\UseTheoremCounterForNextEquation
\begin{equation}\label{GammaAdMInclus3}
\Gamma^{\psi'}_{\ad,N'} \into \Gamma^{\psi'}_\ad \ltimes \Hom_{S'}(N',T_\ad(\psi')).
\end{equation}
Since $(B,K',\psi')$ is primitive, this embedding is open by Theorem \ref{MainPrimitive}. The isomorphisms from (\ref{ChangingAGalois}) 
%Then the diagram (\ref{ChangingAGalois}) for the inclusion $A'\supset B$ becomes 
%\UseTheoremCounterForNextEquation
%\begin{equation}\label{ChangingAGalois1}
%\vcenter{\xymatrix@C-13pt@R-5pt{ 
%\Aut_{A'_\ad}(T_\ad(\phi')) & \Gamma^{\phi'}_\ad \ar@{}[l]|-\supset \ar@{=}[d]^-{\wr} 
%&& \ar@{->>}[ll] \Gamma^{\phi'}_{\ad,M'} \ar@{=}[d]^-{\wr} & \ar@{}[l]|\supset
%\Delta^{\phi'}_{\ad,M'} \ar@{^{ (}->}[rr] \ar@{=}[d]^-{\wr}
%&& \Hom_{A'}(M',T_\ad(\phi')) \ar@{^{ (}->}[d] \\
%\Aut_{B_\ad}(T_\ad(\psi')) & \Gamma^{\psi'}_\ad  \ar@{}[l]|-\supset
%&& \ar@{->>}[ll] \Gamma^{\psi'}_{\ad,M'} & \ar@{}[l]|\supset
%\Delta^{\psi'}_{\ad,M'} \ar@{^{ (}->}[rr] && \Hom_B(M',T_\ad(\psi')) \rlap{.}\\}}
%\end{equation}
thus yield an open embedding
\UseTheoremCounterForNextEquation
\begin{equation}\label{GammaAdMInclus4}
\Gamma^{\phi'}_{\ad,N'} \into \Gamma^{\phi'}_\ad \ltimes \Hom_{S'}(N',T_\ad(\phi')).
\end{equation}
Since $N'=S'M'$, the Galois action on $T_\ad(\phi',N')$ is completely determined by the action on $T_\ad(\phi',M')$; in other words the restriction induces a natural isomorphism $\Gamma^{\phi'}_{\ad,N'} \cong \Gamma^{\phi'}_{\ad,M'}$. This together with Proposition \ref{QuasiHom} (a) yields an open embedding
% we can replace $\Hom_{S'}(N',\ \ )$ by $\Hom_{(S')}(M',\ \ )$.
\UseTheoremCounterForNextEquation
\begin{equation}\label{GammaAdMInclus5}
\Gamma^{\phi'}_{\ad,M'} \into \Gamma^{\phi'}_\ad \ltimes \Hom_{(S')}(M',T_\ad(\phi')).
\end{equation}

The next natural step would be the passage from $(K',\phi',M')$ to $(K',\phi'|A,M')$. However, this runs into the problem that $S'$ does not necessarily act on $T_\ad(\phi'|A)$, because $T_\ad(\phi'|A)$ is obtained from $T_\ad(\phi') \cong \prod_{\Fp'\not=\Fp'_0}T_{\Fp'}(\phi')$ by removing all factors with $\Fp'|\Fp'_0$, which are not necessarily preserved by the non-commutative ring~$S'$. 
%For every prime $\Fp\not=\Fp_0$ of $A$ there is a natural isomorphism $T_\Fp(\phi'|A) \cong \prod_{\Fp'|\Fp}T_{\Fp'}(\phi')$. Thus 
%\UseTheoremCounterForNextEquation
%\begin{equation}\label{Tphiphiprime}
%T_\ad(\phi'|A) \cong \prod_{\Fp'\nmid\Fp_0}T_{\Fp'}(\phi')
%\end{equation}
%is obtained from $T_\ad(\phi') = \prod_{\Fp'\not=\Fp'_0}T_{\Fp'}(\phi')$ by removing the factors corresponding to the finitely many primes $\Fp'\not=\Fp'_0$ of $A'$ that divide~$\Fp_0$.
Thus if $\Fp'_0$ is not the only prime above~$\Fp_0$, it would be ugly to precisely describe the image of $\Hom_{(S')}(M',T_\ad(\phi'))$ in $\Hom_A(M',T_\ad(\phi'|A))$ in general, though of course it can be done. We therefore restrict ourselves to two special cases, with the following results:

%%%%%%%%%%%%%%%%%%%%%%%%%%%%%%%%%%%%%%%%%%%%%%%

\begin{Thm}\label{MainCentral}
Assume that $A$ is the center of $\End_{K^\sep}(\phi)$. With $A'\defeq A$ let $(K',\phi',f,B)$ be as in Proposition \ref{PrimRed} (c), and set $S:= \End_{K^\sep}(\phi|B)$. Let $M$ be a finitely generated torsion free $A$-submodule of~$K$. Then $\Delta_{\ad,M}$ is commensurable with the subgroup $\Hom_{(S)}(M,T_\ad(\phi))$ of $\Hom_A(M,T_\ad(\phi))$, and $\Gamma_{\ad,M}$ is commensurable with $\Gamma_\ad \ltimes \Hom_{(S)}(M,T_\ad(\phi))$.
\end{Thm}

\begin{Proof}
With the above notation $f$ induces an isomorphism $M\stackrel{\sim}{\to} M'$ and an embedding of finite index $T_\ad(\phi)\into T_\ad(\phi')$. It also induces an isomorphism ${S\otimes_B\Quot(B)} \allowbreak \cong S'\otimes_B\Quot(B)$ under which the intersection of $S$ and $S'$ has finite index in both. Since $T_\ad(\phi)$ and $T_\ad(\phi')$ are torsion free $A_\ad$-modules, this implies the equalities and an inclusion of finite index in the diagram
$$\xymatrix@C-15pt{
\ \ \Hom_{(S)}(M,T_\ad(\phi))\ \ar@{=}[r]
& \ \Hom_{(S\cap S')}(M,T_\ad(\phi))\ \ar@{^{ (}->}[d]^{f\circ(\ )\circ f^{-1}} \\
\Hom_{(S')}(M',T_\ad(\phi'))\ \ar@{=}[r]
& \ \Hom_{(S\cap S')}(M',T_\ad(\phi'))\rlap{.}\ \\}$$
With (\ref{DeltaAdMInclus2M'}) and the open embedding (\ref{GammaAdMInclus5}) this shows that the image of $\Gal(K^\sep/K')$ associated to $(K',\phi,M)$ is an open subgroup of $\Gamma_\ad \ltimes \Hom_{(S)}(M,T_\ad(\phi))$. This implies the assertion about $\Gamma_{\ad,M}$, from which the assertion about $\Delta_{\ad,M}$ directly follows.
\end{Proof}

%%%%%%%%%%%%%%%%%%%%%%%%%%%%%%%%%%%%%%%%%%%%%%%

In the other special case we drop all assumptions on endomorphisms, but instead assume something about~$M$:

\begin{Thm}\label{MainGeneralSpec}
Let $(A',K',\phi',f,B)$ be as in Proposition \ref{PrimRed} (c), and set $S:= \End_{K^\sep}(\phi'|B)$. Let $M$ be a finitely generated torsion free $A$-submodule of~$K$ such that the natural map 
$$\textstyle S\otimes_AM\to K^\sep,\ \ 
{\sum_i s_i\otimes m_i} \mapsto \sum_i s_i f(m_i)$$ 
is injective. Then $\Delta_{\ad,M}$ is an open subgroup of $\Hom_A(M,T_\ad(\phi))$, and $\Gamma_{\ad,M}$ is an open subgroup of $\Gamma_\ad \ltimes \Hom_A(M,T_\ad(\phi))$.
\end{Thm}

\begin{Proof}
With the above notation the assumption implies that the natural map 
$$\textstyle S'\otimes_Af(M)\to S'f(M) = N',\ \ 
{\sum_i s'_i\otimes m'_i} \mapsto \sum_i s'_i m'_i$$ 
is injective and therefore an isomorphism. Thus by Proposition \ref{QuasiHom} the restriction induces a natural isomorphism 
$$\Hom_{S'}(N',T_\ad(\phi')) \stackrel{\sim}{\to} \Hom_{A}(f(M),T_\ad(\phi')).$$
The openness of the embedding (\ref{GammaAdMInclus4}), together with the surjectivity in (\ref{ChangingAGalois}) for the inclusion $A\subset A'$, thus implies the openness of the embedding
$$\Gamma^{\phi'|A}_{\ad,f(M)} \into \Gamma^{\phi'|A}_\ad \ltimes \Hom_A(f(M),T_\ad(\phi'|A)).$$
With (\ref{DeltaAdMInclus2M'}) it follows that 
$$\Gamma_{\ad,M} \into \Gamma_\ad \ltimes \Hom_A(M,T_\ad(\phi))$$ 
is an open embedding. This proves the assertion about $\Gamma_{\ad,M}$, from which the assertion about $\Delta_{\ad,M}$ directly follows.
\end{Proof}

%%%%%%%%%%%%%%%%%%%%%%%%%%%%%%%%%%%%%%%%%%%%%%%%%%%%%%%%%%%%%%%%%%%%%%%%%%%%%%%%%%%%%%%%%%%%%

% \newpage

\end{document}